\documentclass[12pt]{amsart}
\usepackage{epsfig}

 \newcommand{\new}{\newcommand}                        
 \new{\trunc}{\hat{\otimes}}                           
 \new{\tnsr}{\otimes}                                  
 \new{\tensor}{\otimes}                                
 \new{\iso}{\cong}                                     
 \new{\union}{\cup}                                    
 \new{\W}{\mathfrak{W}}                                
 \new{\g}{\mathfrak{g}}                                
 \new{\h}{\mathfrak{h}}                                
 \new{\m}{\mathfrak{m}}                                
 \new{\n}{\mathfrak{n}}                                
 \new{\p}{\mathfrak{p}}                                
 \new{\id}{\mathbf{1}}                                 
 \new{\de}{\mathfrak{d}}                               
 \new{\e}{\mathfrak{e}}                                
 \new{\dd}{\hat{\de}}                                  
 \new{\ee}{\hat{\e}}                                   
 \new{\ii}{\mathfrak{i}}                               
 \new{\Vect}{\operatorname{Vect}}                      
 \new{\BC}{\mathfrak{B}\mathfrak{C}}                   
 \new{\uqg}{U_q(\g)}                                   
 \new{\bracket}[1]{\langle#1\rangle}                   
 \new{\qdim}{\operatorname{qdim}}                      
 \new{\Hom}{\operatorname{Hom}}                        
 \new{\C}{\mathbb{C}}                                  
 \new{\FF}{\mathbb{F}}                                 
 \new{\R}{\mathbb{R}}                                  
 \new{\Cat}{\mathfrak{C}}                              
 \new{\ZZ}{Z}                                          
 \new{\F}{\mathfrak{F}}                                
 \new{\Z}{\mathbb{Z}}                                  
 \new{\Q}{\mathbb{Q}}                                  
 \new{\N}{\mathbb{N}}                                  
 \newcounter{letter}                                   
 \newenvironment{alist}{
 \begin{list}{{(\alph{letter})}}{\usecounter{letter}}
 }{\end{list}}                                         
 \newtheorem{thm}{Theorem}                             
 \newtheorem{theorem/definition}{Theorem/Definition}   
 \newtheorem{prop}{Proposition}                        
 \newtheorem{cor}{Corollary}                           
 \newtheorem{lem}{Lemma}                               
\theoremstyle{definition}                              
\newtheorem{rem}{Remark}                               
\newtheorem{definition}{Definition}                    

\new{\pic}[5]{\raisebox{#3pt}{
\hspace{#4pt}\epsfig{file=#1.ps,height=#2pt}\hspace{#5pt}}}

\begin{document}
\title{Three-Dimensional $\mathbf{2}$-Framed TQFTs and Surgery}
\author{Stephen F. Sawin}
\address{Fairfield University\\(203)254-4000x2573\\
ssawin@fair1.fairfield.edu}

\begin{abstract}
    The notion of $2$-framed three-manifolds is defined. The 
    category of $2$-framed cobordisms is described, and used to define 
    a $2$-framed three-dimensional TQFT.   Using skeletonization and 
    special features of this category,  a small set of data and 
    relations is given that suffice to construct a $2$-framed 
    three-dimensional TQFT.   These data and relations are expressed 
    in the language of surgery.
\end{abstract}
\maketitle

\section*{Introduction}
The work of Jones \cite{Jones85}, Witten \cite{Witten89a}, Atiyah
\cite{Atiyah89}  and of Re\-sheti\-khin and Turaev \cite{RT90,RT91}
set the stage for how we have approached the Chern-Simons (or
Jones-Witten or Reshetikhin-Turaev) invariants and their cousins ever
since.  The `natural definition' of the invariants is as the partition
function of a quantum gauge theory with the Chern-Simons functional
as action, computed formally by a mathematically inaccessible path
integral.  Witten argued that these invariants could be effectively
computed by techniques of physics, chief among them the cut-and-paste
properties expected of the partition function of a topological quantum
field theory.  Atiyah translated these cut-and-paste properties into a
precise set of axioms which one could prove  a rigorously defined
invariant satisfies,  thus gaining mathematically respectable access
to some of the path integral arguments that are so famously effective
in the hands of physicists.  The term `topological quantum field theory' or
TQFT is generally used in the mathematical literature to refer to
anything which satisfies these axioms, though of course they capture
only a portion of the tools physicists use to understand TQFTs.
Finally Reshetikhin and Turaev used quantum group techniques to
construct a rigorous topological invariant corresponding to
Chern-Simons theory with gauge group $SU(2)$  and showed that it
satisfied a weaker  version of Atiyah's axioms  (in the language of
this paper, they showed everything but the mending axiom). 

 The
current state of the art is  that from every modular
category one can construct a TQFT  (See Turaev \cite{Turaev94}).  A
variety of interesting invariants fit 
into this framework, most of them  closely related to Chern-Simons
theory and, like Chern-Simons theory, are three-dimensional and
involve an additional topological structure called a \emph{$2$-framing}
or \emph{biframing.}   Many invariants do not fit this framework, but appear to be
analogous in some sense, such as those of Hennings \cite{Hennings96},
Lyubashenko \cite{Lyubashenko95b,Lyubashenko95d} and Kuperberg
\cite{Kuperberg96}  (see also related work of Kerler
\cite{Kerler95,Kerler97,Kerler98b} and Kauffman and Radford \cite{KR95b}).  

Because the notion of TQFT is fairly involved, and the topology and
algebra required to discuss it is complex  (though certainly not deep),
it is often easier to deal with these invariants purely as
three-manifold invariants, ignoring the underlying TQFT structure.
Likewise, it is easier to ignore the technicalities of the $2$-framing
by including a correction term to make the invariant of closed
manifolds independent of the $2$-framing and treating the associated
representation of the $2$-framed mapping class group as a projective
representation of the ordinary mapping class group.  Such
simplifications, however, leave out important
topological information, and perhaps more importantly, sidestep some
of the best connections available between the well worked out topology
and combinatorics of constructing invariants and the deep but very
poorly understood geometry and physics underlying these invariants.

This article reduces the definition of a three-dimensional $2$-framed
TQFT to a small set of relatively straightforward conditions to check,
and then expresses these conditions in the language of surgery (many
of the invariants are computed and defined most naturally in terms of
surgery presentations of $2$-framed three-manifolds).  There are two 
motivations for this.  The first is that these conditions offer 
simplifications and fresh insight in the  construction of TQFTs associated to
the Reshetikhin-Turaev invariants, and more generally constructing TQFTs from
modular categories.  In particular, \cite{Sawin96a} gives a relatively brief
self-contained construction of TQFTs from modular categories using these
techniques. In fact the main results here are used in that paper with a
reference to unpublished lecture notes, and the present article offers a
published reference for them.   The second motivation is as a 
framework for generalizations of the notion of TQFT, such as  adding 
spin and other structures, extending the theory to encompass larger 
codimensions as in Walker \cite{Walker??}, and weakening the TQFT 
structure to apply to the Hennings and Kuperberg invariants.

The first section reviews the notion of a $2$-framing of a
three-manifold, generalizes it to manifolds with boundary, constructs
the symmetric monoidal category of $2$-framed cobordisms, and defines a
TQFT in terms of this category.  The second section simplifies this
definition by general category theoretic arguments and use of
structure special to the $2$-framed cobordism category.  The third
section translates this simplified set of conditions into the language
of surgery, in terms of which many of the interesting invariants are
defined.  

\section{$2$-Framings and the $2$-Framed Cobordism Category}

\subsection{$2$-Framings}
Atiyah \cite{Atiyah90} defined the notion of a $2$-framing  and
analyzed its properties, identifying it as 
the anomaly in Witten's construction of Chern-Simons theory.  Several 
different but equivalent formulations of this anomaly have been given, 
most centering either on  bundles over a three manifold as in Atiyah's 
original definition or on a choice of four-manifold which the  
three-manifold  bounds.  We follow the second approach, because of the 
simplicity of the definition, the underlying four-dimensional nature 
of Chern-Simons theory (see Dijkgraaf and Witten \cite{DW90}), and the 
surgery description.   Our approach follows most closely  Walker \cite{Walker??}.  The chief disadvantage 
is the awkwardness of viewing a four-dimensional category as a 
three-dimensional category with extra structure.

All of our manifolds will be assumed to be compact, oriented, and 
smooth, but 
not necessarily connected or closed (although we will usually add the term 
``with boundary'' when speaking of manifolds which may not be closed).

A \emph{$2$-framed three-manifold} is a four-manifold $M$ with 
boundary, with the proviso that each component of $M$ have at most one 
component of boundary (in fact since we will only really be concerned 
with our four-manifolds up to cobordism, this does not affect the 
theory, but makes a number of technicalities easier to deal with).   
We will refer to  $M$ as a choice of $2$-framing on its boundary 
$\partial M.$  A \emph{$2$-framed diffeomorphism} between $2$-framed 
three-manifolds $M$ and $N$ is a diffeomorphism $f\: \partial M \to 
\partial N$ together with a five-manifold $W$ whose boundary is $M 
\cup_f N,$  the closed manifold formed by identifying $\partial M$ 
with $\partial N$ via $f.$   Thus \cite{Kirby89} two $2$-framed three-manifolds 
are $2$-framed diffeomorphic if and only 
if their boundaries are diffeomorphic and they have the same signature.

A \emph{$2$-framed surface $\Sigma$} is a three-manifold with 
boundary (again with the proviso that each component has connected 
boundary) , and 
we will refer to $\Sigma$ as a choice of $2$-framing on $\partial 
\Sigma.$   A \emph{$2$-framed diffeomorphism}  of 
$2$-framed surfaces is just a diffeomorphism of the $2$framed 
surfaces as three-manifolds.     A \emph{$2$-framed three-manifold with 
boundary $(M,\Sigma)$} is 
a $2$-framed three-manifold $M$ together with a $2$-framed surface 
$\Sigma$ which as a three-manifold is a submanifold of the boundary $\partial 
M.$  We will refer to $(M,\Sigma)$ as a choice of $2$-framing on 
$\partial M-\Sigma.$   A 
\emph{diffeomorphism of $2$-framed manifolds} $(M,\Sigma) $
and $(N,\Gamma)$ with boundary  is a diffeomorphism of the $2$-framed 
three manifolds whose underlying ordinary diffeomorphism of the 
boundaries sends $\Sigma$ to $\Gamma.$  

\begin{figure}[hbt]
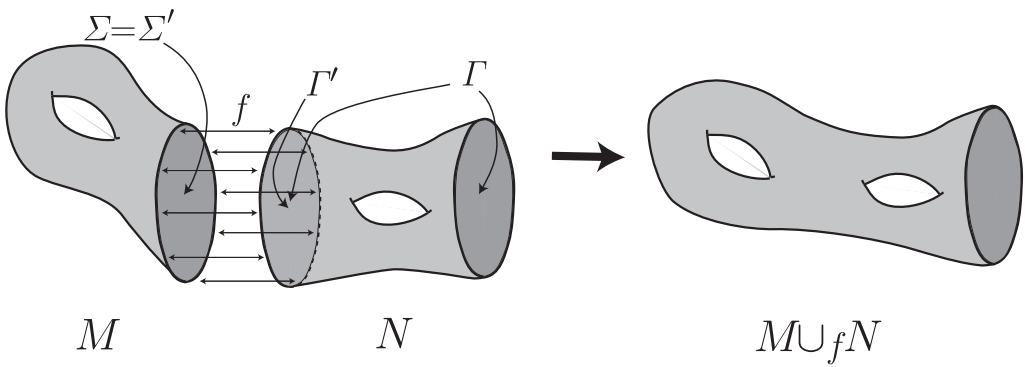

    $$\pic{gluing}{135}{-75}{-15}{0}$$
    \caption{The gluing of $M$ and $N$} \label{fg:gluing}
\end{figure}

If $(M,\Sigma)$ and $(N,\Gamma)$ are $2$-framed manifolds with 
boundary and $f$ is a $2$-framed diffeomorphism of subsets of their boundaries $\Sigma' 
\subset \Sigma$ 
and $\Gamma' \subset \Gamma',$ then we can we can form the gluing of 
these, $M \cup_f N,$  which is the $2$-framed three-manifold with 
boundary which as a  four manifold is $M$ and $N$ glued together along the 
subsets $\Sigma'$ and $\Gamma'$ of their boundary, and whose boundary 
$2$-framed  surface is the subset $\Sigma-\Sigma' \cup 
\Gamma-\Gamma'.$  See Figure \ref{fg:gluing} for a pictorial 
representation of gluing (here for visual clarity all dimensions have 
been reduced by one, so that $M$ and $N$ which should be represented 
as four-manifolds are pictured as three-manifolds, $\Sigma$ and 
$\Gamma$ which should be represented by three-manifolds are pictured 
as two-manifolds, etc.   

Notice if $(M,\Sigma)$ and $(N,\Gamma)$ are 
$2$-framed three-manifolds with boundary, then $\Sigma$ is a choice 
of $2$-framing on the boundary of $\partial M,$ $\Gamma$ is a choice 
of $2$-framing on the boundary of $\partial N,$ and $M \cup_f N$ is a 
choice of $2$-framing on $(\partial M -\Sigma)  \cup_{\partial f} (\partial 
N -\Gamma).$  Thus 
roughly speaking the notion of  ``a $2$-framing on'' commutes with the 
operations of taking the boundary and gluing.

Below we shall refer to various notions for $2$-framed manifolds and 
diffeomorphisms, 
such as disjoint union and 
orientation reversal,   which are defined exactly as for ordinary 
manifolds and diffeomorphisms, and we will not bother to state the definition 
precisely.

\subsection{The $2$-framed cobordism category}
Let us define a $2$-framed cobordism as $\m=(M,\Sigma,\Gamma,f),$
where $M$ is a $2$-framed three-manifold with boundary, $\Sigma$ and
$\Gamma$ are $2$-framed two-manifolds, and $f$ is an orientation preserving
$2$-framed diffeomorphism of $\Sigma^* \cup \Gamma$ to the boundary 
of $M$ (here
$\Sigma^*$ means $\Sigma$ equipped with the opposite orientation).   We
say that two $2$-framed cobordisms $\m=(M,\Sigma,\Gamma,f)$ and
$\m'=(M', \Sigma, \Gamma, f')$ are the same if there is a
$2$-framed diffeomorphism $F\colon M \to M'$ such that $F \circ 
f=f'$ (in particular, they are the same if and only if they are 
$2$-framings of diffeomorphic three-dimensional ordinary cobordisms 
and their underlying four-manifolds have the same signature). If
$\m=(M,\Sigma, \Gamma, f)$ we shall simply write $\m\colon\Sigma \to \Gamma$ in
the frequent situation that we do not need to refer explicitly to the underlying
manifold and parameterization map $f.$

If $\m=(M,\Sigma,\Gamma,f)$ and $\n=(N,\Gamma,\Delta,g)$ are
two $2$-framed cobordisms, then $h=g|_{\Gamma} \circ (f|_\Gamma)^{-1}$
is an orientation-reversing $2$-framed diffeomorphism between subsets of
the boundary of $M$ and $N,$ and thus we get a $2$-framed three-manifold
$M \cup_h N$ and a $2$-framed cobordism, the \emph{composition} of $\m$
and $\n$ defined as 
$$\n \circ \m = (M \cup_h N, \Sigma, \Delta, f|_{\Sigma^*} \cup
g|_{\Delta}).$$

If $\Sigma$ is a $2$-framed two-manifold, $(\Sigma \times I, \Sigma^* 
\times \{0\} \cup \Sigma \times \{1\})$ is of 
course a $2$-framed three-manifold with boundary.
This determines a $2$-framed cobordism
$$\id_\Sigma=(\Sigma \times I, \Sigma,\Sigma, f)$$
where $f$ sends the first $\Sigma$ to $\Sigma \times \{0\}$ and the
second to $\Sigma \times \{1\},$ both by the identity map.  One
readily sees that 
$$(\m \circ \n) \circ \p= \m \circ (\n \circ \p)$$
$$\id_\Sigma \circ \m= \m \circ \id_\Sigma=\m$$
whenever the compositions are well-defined.  

If we  take the somewhat
pedantic position that the empty three-manifold $\emptyset$ is a
$2$-framed two-manifold and the empty map $\emptyset$ is a $2$-framed diffeomorphism
from $\emptyset$ to $\emptyset,$ we can view closed $2$-framed
three-manifolds $M$ as cobordisms via
$(M,\emptyset,\emptyset,\emptyset).$  For a slight additional cost in
pedantry, we may also view $\emptyset$ as a $2$-framed three-manifold,
allowing us to shamelessly write
$\id_\emptyset=(\emptyset,\emptyset,\emptyset,\emptyset).$  The
relations above still hold for these vacuous cobordisms.

If $\m=(M,\Sigma_1,\Sigma_2,f)$ and $\n=(N,\Gamma_1,\Gamma_2,g)$ are
two $2$-framed cobordisms,  the disjoint union is
$$\m \cup \n=(M \cup N, \Sigma_1 \cup \Gamma_1,\Sigma_2 \cup
\Gamma_2,f \cup g).$$
Again it is clear that
$$\id_\emptyset \cup \m= \m \cup \id_\emptyset=\m,$$
$$(\m \cup \n) \cup \p=\m \cup(\n \cup \p)\,\,\,\text{and}$$
$$(\m \cup \n) \circ (\m' \cup \n')= (\m \circ \m') \cup (\n \circ
\n').$$

If $\Sigma$ and $\Gamma$ are any two $2$-framed two-manifolds define
$$\sigma_{\Sigma,\Gamma}=((\Sigma \cup \Gamma) \times I, \Sigma \cup
\Gamma, \Gamma \cup \Sigma, f)$$
where $f$ maps $\Sigma \cup \Gamma$ to $(\Sigma \cup \Gamma) \times
\{0\}$ by the identity map and and $\Gamma \cup \Sigma$ to $(\Sigma
\cup \Gamma) \times \{1\}$ by the obvious flip of the identity.  Of
course
$$\sigma_{\Sigma,\Gamma} \circ (\m \cup \n) \circ
\sigma_{\Gamma,\Sigma}=\n \cup \m\,\,\,\text{and}$$
$$(\id_\Gamma \cup \sigma_{\Sigma,\Delta}) \circ (\sigma_{\Sigma,\Gamma}
\cup \id_{\Delta})=\sigma_{\Sigma,\Gamma \cup \Delta}.$$

This is all to say that there is a category whose objects are $2$-framed
two-manifolds and whose  morphisms are $2$-framed cobordisms with
composition defined above and $\id_\Sigma$ as identity, that $\cup$
forms a strict monoidal product on this category with $\emptyset$ as
trivial object, and that $\sigma$ forms a symmetry for this monoidal
category.  We call this symmetric monoidal category the \emph{$2$-framed
cobordism category} $\BC.$

$\BC$ has some particularly nice properties, revolving around the
artificiality of dividing up the boundary into two pieces.  In
particular to each $2$-framed two-manifold $\Sigma$ we can associate
$\Sigma^*,$ which is $\Sigma$ with the opposite orientation, and cobordisms
$$\de_\Sigma =(\Sigma \times I, \Sigma \cup \Sigma^*, \emptyset,f)
\,\,\,\text{and}$$
$$\e_\Sigma=(\Sigma \times I,\emptyset, \Sigma^* \cup \Sigma, f)$$
where in both cases $f$ is the obvious map from $\Sigma^* \cup \Sigma$
to the boundary of $\Sigma \times I.$  These morphisms have the
properties
$$(\de_\Sigma \cup \id_\Sigma) \circ (\id_\Sigma \cup
\e_\Sigma)=\id_\Sigma\,\,\,\text{and}$$
$$(\id_{\Sigma^*} \cup \de_\Sigma) \circ (\e_\Sigma \cup
\id_{\Sigma^*})=\id_{\Sigma^*}.$$ 
A monoidal category which admits such morphisms is called a \emph{rigid
monoidal category.}  A rigid symmetric monoidal category is sometimes
referred to in the literature as a closed category or a tensor
category, though the latter term
is often used as a synonym for monoidal category.

\begin{rem}
The reader may note the asymmetry of the treatment of the monoidal and
symmetric structure on the one hand, which are part of the definition of $\BC,$ and the
rigidity on the other, which is observed after the fact to be a structure that can
be assigned to the category.  When we construct functors from this
category to other categories, the target categories will have explicit
symmetric monoidal structures, and we will want the functors to make these
structures correspond exactly, but will not ask this of the rigidity.
The reason is that while the symmetry and rigidity are both quite
canonical in this category, in the linear category to which we will
map, the rigidity is distinctly fussier than the symmetry.
\end{rem}

\subsection{$2$-Framed TQFTs}
The category $\Vect$ whose objects are finite-dimensional vector spaces
and whose morphisms are linear maps with the usual composition can
also be made into a symmetric monoidal category.   The monoidal
structure sends vector spaces $V$ and $W$ to the tensor product $V
\tensor W$ and linear maps $f\colon V \to V'$ and $g\colon W \to W'$ to the
tensor product $f \tensor g\colon V \tensor W \to V' \tensor W',$ and makes the
ground field  the trivial object.  The symmetric morphism
$\sigma_{v,w}$ sends the vector $v \tensor w$ to $w \tensor v.$

A symmetric monoidal functor between two symmetric monoidal categories
is an ordinary functor which preserves the symmetry and monoidal
structure.
\begin{rem}
A more careful treatment is to assume both categories are \emph{weak}
monoidal categories.  For instance, instead of asserting that vector
spaces $(X \tensor Y) \tensor Z$ and $X \tensor (Y \tensor Z)$ are
equal, we should provide a well-behaved isomorphism between them.  Any
careful attempt to define either union of manifolds or tensor product
of vector spaces requires this.   Nevertheless, it is the common and
generally harmless custom in mathematics to view the vector spaces
$(X \tensor Y) \tensor Z$ and $X \tensor (Y \tensor Z)$ as ``the
same,'' and in the interest of clarity and simplicity we will continue
this relaxed practice, pointing out occasional subtleties.
\end{rem}

\begin{definition}
A TQFT is a symmetric monoidal functor from $\BC$ to $\Vect.$
\end{definition}

Specifically, a TQFT $\ZZ$ assigns a vector space $\ZZ(\Sigma)$ to
each $2$-framed two-manifold $\Sigma,$ and a linear map
$\ZZ(\m)\colon\ZZ(\Sigma) \to \ZZ(\Gamma)$ to each $2$-framed cobordism
$\m\colon\Sigma \to \Gamma,$ satisfying the following conditions:
\begin{enumerate}
\item $\ZZ(\m) \circ \ZZ(\n)= \ZZ(\m \circ \n)$
\item $\ZZ(\id_\Sigma)=\id_{\ZZ(\Sigma)}$
\item $\ZZ(\Sigma \cup \Gamma) = \ZZ(\Sigma) \tensor \ZZ(\Gamma)$
\item $\ZZ(\m \cup \n)=\ZZ(\m) \tensor \ZZ(\n)$
\item $\ZZ(\emptyset)= \FF,$ the ground field
\item $\ZZ(\sigma_{\Sigma,\Gamma})=\sigma_{\ZZ(\Sigma),\ZZ(\Gamma)}.$
\end{enumerate}

Two TQFTs are \emph{equivalent} if there is a natural isomorphism
between them.  That is, for each $2$-framed surface $\Sigma,$ an
isomorphism $i_\Sigma\colon\ZZ_1(\Sigma) \to \ZZ_2(\Sigma)$ such that if
$\m\colon\Sigma \to \Gamma$ then $\ZZ_2(\m)\circ i_\Sigma=i_\Gamma \circ
\ZZ_1(\m)$ and $i_{\Sigma \cup \Gamma}=i_\Sigma \tensor i_\Gamma.$

\begin{rem}
    Notice that a TQFT offers in particular an invariant of cobordisms 
    from $\emptyset$ to $\emptyset,$ which is to say an invariant of 
    closed $2$-framed three-manifolds up to $2$-framed diffeomorphism. 
    This means   an invariant of four-manifolds with 
    boundary up to cobordism.  
\end{rem}

\section{Characterization of TQFTs}

\subsection{Skeletonization}
One problem with defining a TQFT is that we seem to have a huge array
of objects, each of which would need to be assigned a vector space
simply as a starting point.  Of course it is clear that many of these
are not fundamentally different, and there is much less significant
choice involved then it would appear.  This is a common situation in
category theory, where the notion of isomorphism captures it
precisely.  Specifically, two objects $\Sigma$ and $\Gamma$ are
\emph{isomorphic} if there exist morphisms $\m\colon\Sigma \to \Gamma$ and $\n\colon\Gamma \to \Sigma$
such that $\m \circ \n=\id_\Gamma$ and $\n \circ \m = \id_\Sigma.$

\begin{lem}\label{lm:skeleton}
Two $2$-framed surfaces are isomorphic as objects in $\BC$ if and only 
if they are $2$-framings of diffeomorphic $2$-manifolds.
\end{lem}

\begin{proof}
    Of course,  an isomorphism between two $2$-framed surfaces must 
    be a $2$-framing of an isomorphism between two ordinary surfaces 
    in the category of ordinary three-dimensional cobordisms, and thus 
    the surfaces must be diffeomorphic.
    
Suppose $\Sigma $ and $\Gamma$ are $2$-framed surfaces  which are 
$2$-framings of the same surface $S$  (here we have absorbed the 
ordinary diffeomorphism of the surfaces for notational simplicity).  
Viewing $\Sigma$ and $\Gamma$ as three-manifolds with boundary $S,$ 
we form the closed three-manifold $\Sigma^* \cup_S (S \times I) 
\cup_S \Gamma$ and choose some four-manifold $W$ which it bounds.  
Then $(W, \Sigma^* \cup \Gamma)$ is a $2$-framed three-manifold with 
boundary and $\n=(W, \Sigma, \Gamma, f)$ is a cobordism from $\Sigma$ 
to $\Gamma$ with the obvious $f.$  Likewise glue $\Gamma^* \cup_S (S 
\times I) \cup_S \Sigma$ and choose a $2$-framing to get a $2$-framed 
cobordism $\n' \colon \Gamma \to \Sigma.$   Now $\n' \circ \n$  
represents a 
$2$-framing on $S\times I,$ and thus 
is the same $2$-framed cobordism as $\id_\Sigma$ if and only if it has 
the same signature (here by the signature we mean the signature of the 
four-manifold which forms the $2$-framed manifold underlying the 
$2$-framed cobordism).  Now by Wall's nonadditivity of the signature 
result \cite{Wall69} the signature of $\n' \circ \n$ is not the sum of the 
signatures of $\n$ and $\n',$ but the difference depends only on the 
boundaries, so by connect summing an appropriate number of  copies of 
$\pm CP^2$ to $\n'$ we get a new cobordism $\n^{-1}$ such that 
$\n^{-1} \circ \n$ has the same signature as, and thus is the same 
morphism as, $\id_\Sigma.$  Therefore $\Sigma$ and $\Gamma$ are 
isomorphic.
\end{proof}

Thus there is one isomorphism class of objects for each genus.

Recall that two (symmetric, monoidal) categories $\Cat$ and $\Cat'$
are \emph{equivalent} if there are a pair of (symmetric, monoidal)
functors $F,$ $G$ between them such that $F\circ G$ and $G \circ F$
are naturally isomorphic to the identity.

\begin{prop}
For any choice of a $2$-framed surface for each genus, $\BC$ is
equivalent to the full subcategory $\BC'$ formed from arbitrary unions
of these surfaces  (i.e., arbitrary unions of these surfaces are the
objects, all morphisms between 
them form the morphisms).
\end{prop}

\begin{proof}
The functor from $\BC'$ to $\BC$ is just the inclusion functor.  For
the other direction choose for each object of $\BC$ which is
connected an isomorphism to the chosen object of the same genus (using
Lemma \ref{lm:skeleton}).  An arbitrary object $\Sigma$ of $\BC$ is a union
of connected pieces, so it is the domain of  an isomorphism $\ii_\Sigma$
which is a union of the chosen isomorphisms.  Then the functor sends
$\Sigma$ to the isomorphic object in $\BC'$ and a morphism $\m\colon\Sigma \to \Gamma$ to
$\ii_\Gamma \circ \m \circ \ii_\Sigma^{-1}.$  The composition of functors one
way is the identity, and the other way is isomorphic to the identity
via the natural isomorphism $\ii_\Sigma^{-1}.$
\end{proof}

\begin{rem}
Typically the so-called skeletonization process of the previous
proposition results in a weak monoidal structure.  In this case the
monoidal structure is free, in the sense that two-manifolds and
cobordisms can be uniquely written as the union of connected components
(up to symmetry) and because of this there is no weakening of the structure.
\end{rem}

\begin{cor}
Every symmetric monoidal functor  from the full subcategory $\BC'$
given above to $\Vect$ uniquely determines a TQFT up to equivalence.
\end{cor}

\begin{proof}
Let $F$ be the equivalence functor $\BC \to \BC',$ and $\ZZ_0$
be the functor $\BC' \to \Vect.$  Then of course $\ZZ_0 \circ F$ is a
symmetric monoidal functor $\BC \to \Vect,$ and hence a TQFT.  If
$\ZZ$ is a TQFT whose restriction to $\BC'$ is $\ZZ_0,$ and 
 $\ii_\Sigma$ is the natural isomorphism from $F$ to the identity
functor on $\BC,$ then $\ZZ(\ii_\Sigma)$ is an isomorphism from
 $\ZZ_0\circ F(X)$ to $\ZZ(X).$  It is natural because $\ii_\Sigma$ is
natural, and thus $\ZZ$ and $\ZZ_0\circ F$ are equivalent.  
\end{proof}

\subsection{Minimal data and relations}

There is still quite a bit of redundancy involved in checking that
something is a TQFT.  In particular, the division of the parameterized
boundary into source and target seems quite arbitrary:  All such
divisions are the same in some sense, and it seems that once you know
the value of the TQFT on one of these you know them all.  This
subsection makes that notion precise.

Consider the set of all \emph{basic} cobordisms: i.e. $2$-framed cobordisms
$\m\colon \bigcup_{i=1}^n \Sigma_{g_i} \to \emptyset$ with target the
trivial surface, source a union of the chosen $2$-framed surfaces, and
whose underlying four-manifold $M$ is connected.

Choose for each chosen $\Sigma_g$ an orientation-reversing  map
$S\colon\Sigma_g \to \Sigma_g$ such that $S^2=\operatorname{id}$ (the choice of
identification of $\Sigma_g$ with its orientation reversal is necessary in describing
the rigidity structure because we have skeletonized:  Making the
identification an involution ($S^2=\operatorname{id}$) is not, but allows us to keep
track of fewer issues of ordering).  Choose for each $\Sigma_g$
duality morphisms  $\de_g=(\Sigma_g \times I, \Sigma_g \cup \Sigma_g^*,
\emptyset, f)$ and $\e_g=(\Sigma_g \times I, \emptyset, \Sigma_g^*
\cup \Sigma_g,f)$ and compose the parameterization maps on the
$\Sigma_g^*$ pieces with $S$ to get a \emph{cap} morphism $\dd_g$ and a
\emph{cup} morphism $\ee_g$
\begin{align*}
\dd_g\colon&\Sigma_g \cup \Sigma_g \to \emptyset\\
\ee_g\colon&\emptyset \to \Sigma_g \cup \Sigma_g
\end{align*}
such that 
$$((\dd_g \cup \id_{g})\circ (\id_{g} \cup \ee_g)=
(\id_{g} \cup \dd_g)\circ (\ee_g \cup
\id_{g})=\id_{g}$$
where $\id_g$ is shorthand for $\id_{\Sigma_g}.$
Notice that $\dd_g$ is a basic cobordism.

\begin{figure}[hbt]
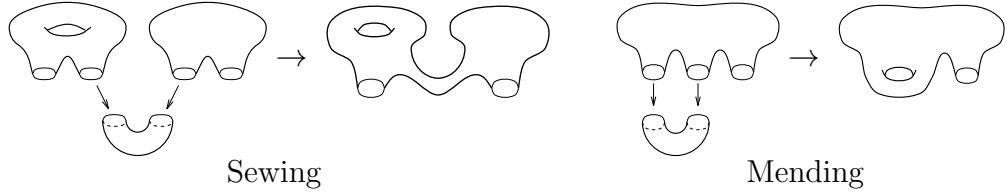

\begin{tabular*}{5.4in}{c@{\extracolsep{\fill}}c}
$\pic{lsewex}{58}{-34}{0}{0} \rightarrow \pic{rsewex}{35}{-12}{0}{0}$
& $\pic{lmendex}{58}{-34}{0}{0} \rightarrow
\pic{rmendex}{35}{-12}{0}{0} $\\
Sewing & Mending
\end{tabular*}
    \caption{Mending versus sewing for ordinary two-manifolds with 
    boundary} \label{fg:mands}
    
    \end{figure}

We wish to consider three operations on basic cobordisms
\begin{itemize}
\item[\textbf{Permuting:}] Suppose $\sigma$ is a permutation of $(1,
\ldots,n).$  If $f\colon\bigcup_{i=1}^n \Sigma_{g_i} \to \partial M$
then $\sigma$ acts on the left in an obvious way on $f,$
$\sigma(f)(x_1,\ldots,x_n)=f(x_{\sigma^{-1}(1)},\ldots,x_{\sigma^{-1}(n)}),$
and thus if $\m=(M,\bigcup_{i=1}^n \Sigma_{g_i},\emptyset,f)$ is
basic we get a new basic cobordism $\sigma(\m)=(M,\bigcup_{i=1}^n
\Sigma_{g_{\sigma^{-1}(i)}}, \emptyset,\sigma(f))$ which of course is
$\m$ composed with an appropriate product of symmetry morphisms.
\item[\textbf{Sewing:}]  Suppose $\m\colon\bigcup_{i=1}^m \Sigma_{g_i}\to
\emptyset$ and $\n\colon\bigcup_{i=1}^n \Sigma_{g'_i}\to
\emptyset$ with $g_m=g'_i.$  Then we can form the \emph{sewing} of
$\m$ and $\n$ 
$$\m \cup_s \n=(\m \cup \n) \circ(\id_{g_1}
\cup \cdots \cup \id_{g_{m-1}} \cup \ee_{g_m} \cup \id_{g'_2} \cup
\cdots \cup \id_{g'_n}).$$
This corresponds to gluing the underlying manifolds together via $S$ along the
$\Sigma_{g_m}$ component of the boundaries.   See Figure \ref{fg:mands} 
for an illustration.
\item[\textbf{Mending}]  Suppose $\m\colon\bigcup_{i=1}^m \Sigma_{g_i} \to
\emptyset$ and $g_1=g_2.$  Then we can form the \emph{mending} of $\m$
$$\m_m=\m\circ(\ee_{g_1} \cup \id_{g_3} \cup \cdots \cup \id_{g_m}).$$
This glues the underlying manifold to itself via $S,$ thus increasing the
dimension of the first homology.  See Figure \ref{fg:mands} for an 
illustration.
\end{itemize}

\begin{thm}\label{th:TQFT}
Suppose we are given for each chosen surface $\Sigma_g$  a
finite-dimensional vector space $\ZZ(\Sigma_g),$ and for each basic
morphism $$\m \colon \bigcup_{i=1}^n \Sigma_{g_i} \to \emptyset$$  a
functional $\ZZ(\m)$ on $\bigotimes_{i=1}^n \ZZ(\Sigma_{g_i})$
satisfying the following conditions (If $\m \colon \emptyset \to \emptyset$
we understand $\ZZ(\m)$ is a functional on the ground field $\FF$).
\begin{alist}
\item \textbf{Symmetry:}  $\ZZ(\sigma(\m))= \sigma(\ZZ(\m)),$ where
$\sigma$ acts on $\bigotimes_{i=1}^n \ZZ(\Sigma_{g_i})$ by the map
$v_1 \tensor \cdots \tensor v_n \mapsto v_{\sigma(1)} \tensor
\cdots \tensor v_{\sigma(n)}$ and hence  on its dual space by
$\sigma(\ZZ(\m))= \ZZ(\m)\circ \sigma^{-1}.$
\item \textbf{Nondegeneracy:}  $\ZZ(\dd_g)$ is a nondegenerate pairing
on $\ZZ(\Sigma_g),$ which we'll call $\bracket{\, ,\,}.$  By
symmetry it is a symmetric pairing, and let us define $\ZZ(\ee_g)$ to
be the canonical dual element $\sum_j v_j \tensor w_j \in
\ZZ(\Sigma_g) \tensor \ZZ(\Sigma_g)$ such that $\sum_j v_j
\bracket{w_j,x}=\sum_j \bracket{x,v_j} \tensor w_j=x.$
\item \textbf{Sewing:} If $\m\colon\bigcup_{i=1}^m \Sigma_{g_i}$ and
$\n\colon\bigcup_{j=1}^n \Sigma_{g'_j}$ with $g_m=g'_1$ then 
$\ZZ(\m) \tensor \ZZ(\n) \circ \Phi=\ZZ(\m
\cup_s \n),$ where 
$$\Phi\colon\bigotimes_{i=1}^{m-1}\ZZ(\Sigma_{g_i}) \tensor
\bigotimes_{j=2}^n \ZZ(\Sigma_{g_j}) \to
\bigotimes_{i=1}^{m}\ZZ(\Sigma_{g_i}) \tensor 
\bigotimes_{j=1}^n\ZZ(\Sigma_{g_j})$$ is the canonical map which is
$\ZZ(\ee_{g_m})$ as a map $\FF \to
\ZZ(\Sigma_{g_m}) \tensor \ZZ(\Sigma_{g_m})$ tensored with the
identity on the other factors.
\item \textbf{Mending:}  If $\m\colon\bigcup_{i=1}^n \Sigma_{g_i}$ with
$g_1=g_2$ then $\ZZ(\m_m)=\ZZ(\m) \circ \Phi$ where
$$\bigotimes_{i=3}^{n} \ZZ(\Sigma_{g_i}) \to
\bigotimes_{i=1}^{n} \ZZ(\Sigma_{g_i}) $$
is the canonical map which is $\ZZ(\ee_{g_1})$ tensored with the
identity on the other factors.
\end{alist}
Then $\ZZ$ determines a TQFT, unique up to equivalence, whose value on
the chosen surfaces and basic morphisms is as given.
\end{thm}

\begin{proof}
If $\Sigma=\bigcup_{i=1}^n \Sigma_{g_i},$ define
$\ZZ(\Sigma)=\bigotimes_{i=1}^n \ZZ(\Sigma_{g_i}).$

For each permutation morphism $\sigma\colon\bigcup_{i=1}^n \Sigma_{g_i} \to
\bigcup_{i=1}^n \Sigma_{g_{\sigma(i)}},$ define
$\ZZ(\sigma)\colon\bigotimes_{i=1}^n \ZZ(\Sigma_{g_i}) \to
\bigotimes_{i=1}^n \ZZ(\Sigma_{g_{\sigma(i)}})$ to be the
corresponding permutation map.  Because the permutation groupoid
embeds into $\BC'$ we have $\ZZ(\sigma_1) \circ
\ZZ(\sigma_2)=\ZZ(\sigma_1 \circ \sigma_2)$ and because of the
symmetry axiom we have $\ZZ(\m) \ZZ(\sigma)=\ZZ(\m \circ \sigma)$ when
$\m$ is basic.

If $\m\colon\bigcup_{i=1}^m \Sigma_{g_i} \to \emptyset$ is not basic because it 
has an underlying manifold which is not connected, then there exists a
permutation $\sigma$ such that $\m \circ \sigma=\m_1 \cup \cdots \cup
\m_n$ is a union of basic morphisms.  Define $\ZZ(\m)=(\ZZ(\m_1)
\tensor \cdots \tensor \ZZ(\m_n)) \circ \ZZ(\sigma^{-1}).$  The choice
of $\sigma$ was not unique, since we could have reordered the factors
$\m_i,$ but it is easy to check that   $\ZZ(\m)$
is well-defined, and that $\ZZ$ on this extended collection of
morphisms satisfies all the axioms in the statement of the theorem as
well as $\ZZ(\m \cup \n)=\ZZ(\m) \tensor \ZZ(\n)$ and $\ZZ(\m \circ
\sigma)= \ZZ(\m) \circ \ZZ(\sigma).$

Now if $\Sigma=\bigcup_{i=1}^n \Sigma_{g_i},$ define 
$\Sigma^\dagger=\bigcup_{i=1}^n \Sigma_{g_{n-i +1}}$ and 
$$\dd_\Sigma=(\Sigma \times I,\Sigma^\dagger \cup \Sigma,\emptyset,f)$$
where  $f$ sends
$\Sigma$ to $\Sigma \times \{0\}$ via the identity map and
$\Sigma^\dagger$ to $\Sigma \times \{1\}$ via the obvious permutation
composed with  the union of a copy of the $S$ map for each
$\Sigma_{g_i}.$  It is easy to check that $\ZZ(\dd_\Sigma)$ is the
pairing on $\ZZ(\Sigma)$ determined by the pairings $\dd_i$ on each
$\Sigma_{g_i}.$  Define $\ee_\Sigma$ so that $(\dd_\Sigma \tensor 
\id_{\Sigma})(\id_\Sigma \tensor \e_\Sigma)=\id_\Sigma$ and 
$(\id_{\Sigma^\dagger} \tensor \dd_\Sigma)(\ee_\Sigma \tensor 
\id_{\Sigma^\dagger})=\id_{\Sigma^\dagger}$ and define
$\ZZ(\ee_\Sigma)$ to be the element of
$\ZZ(\Sigma) \tensor \ZZ(\Sigma^\dagger)$ dual to the pairing.

If $\m\colon\Gamma \cup \Sigma \to \emptyset$ and $\n\colon \Sigma^\dagger
\cup \Delta \to \emptyset$ we can construct $\m \cup_\Sigma \n \colon \Gamma
\cup \Delta \to \emptyset, $ which is $(\m \cup \n) \circ
(\id_{\Gamma} \cup \ee_\Sigma \cup \id_\Delta),$  by a sequence of sewings,
mendings, and permutations. It then follows from the sewing, mending 
and symmetry axioms that
$$\ZZ(\m \cup_\Sigma \n)=(\ZZ(\m) \tensor \ZZ(\n)) \circ \Phi$$
where $\Phi$ is formed by tensoring the identity map with the map
associated to $\ee_\Sigma.$

Finally, if $\m\colon\Sigma \to \Gamma$ is an arbitrary element of $\BC',$
define $\hat{\m} \colon \Gamma^\dagger \cup \Sigma \to \emptyset$ by 
$$\hat{\m}=\dd_\Gamma \circ(\id_{\Gamma^\dagger} \cup \m)$$
noting that
$$\m=(\id_\Gamma \cup \hat{\m}) \circ (\ee_\Gamma \cup \id_\Sigma)$$
\and define
$$\ZZ(\m)=(\id_{\ZZ(\Gamma)} \tensor \ZZ(\hat{\m})) \circ
(\ZZ(\ee_\Gamma) \cup \id_{\ZZ(\Sigma)}).$$ 

We have
only to check that the conditions after the definition of a TQFT are
satisfied.
\begin{enumerate}
\item If $\n\colon\Sigma \to \Gamma$ and $\m\colon\Gamma \to \Delta$
then 
$\ZZ(\m \circ \n)= \ZZ(\m) \circ \ZZ(\n).$  To see this note
\begin{multline*}
\ZZ(\widehat{\m \circ \n})= \ZZ(\hat{\m} \cup_\Gamma \hat{\n})\\
= (\ZZ(\hat{\m}) \tensor \ZZ(\hat{\n})) \circ
(\id_{\ZZ(\Delta^\dagger)} \tensor \ZZ(\ee_\Gamma) \tensor \id_{\ZZ(\Sigma)})\\
=(\ZZ(\dd_\Delta) \tensor \ZZ(\dd_\Gamma)) \circ
(\id_{\ZZ(\Delta^\dagger)} \tensor \ZZ(\m) \tensor
\id_{\ZZ(\Gamma^\dagger)} \tensor \ZZ(\n) ) \circ (\id_{\ZZ(\Delta)} \tensor
\ZZ(\ee_\Gamma) \tensor \id_{\ZZ(\Sigma)})\\
=\ZZ(\dd_{\Delta}) \circ (\id_{\ZZ(\Delta^\dagger)} \tensor (\ZZ(\m) \circ
\ZZ(\n)))
\end{multline*}
and thus
\begin{multline*}
\ZZ(\m \circ \n)=(\id_{\ZZ(\Delta)} \tensor \ZZ(\widehat{\m \circ
\n})) \circ (\ZZ(\ee_\Delta) \tensor \id_{\ZZ(\Sigma)})\\
=\ZZ(\m) \circ \ZZ(\n).
\end{multline*}
\item $\ZZ(\id_\Sigma)=\id_{\ZZ(\Sigma)}.$  By definition
\begin{multline*}
\ZZ(\id_\Sigma)=(\id_{\ZZ(\Sigma)} \tensor
\ZZ(\widehat{\id_\Sigma})) \circ (\ZZ(\ee_\Sigma) \tensor
\id_{\ZZ(\Sigma)})\\
=(\id_{\ZZ(\Sigma)} \tensor
\ZZ(\dd_\Sigma)) \circ (\ZZ(\ee_\Sigma) \tensor
\id_{\ZZ(\Sigma)})\\
=\id_{\ZZ(\Sigma)}.
\end{multline*}
\item $\ZZ(\Sigma \cup \Gamma)= \ZZ(\Sigma) \tensor \ZZ(\Gamma).$
This is by definition.
\item $\ZZ(\m \cup \n)= \ZZ(\m) \tensor \ZZ(\n).$  Suppose that
$\m\colon\Sigma_1 \to \Gamma_1$ and $\n\colon\Sigma_2 \to \Gamma_2,$ so that
$$
\widehat{\m \cup \n}= (\hat{\m} \cup \hat{\n}) \circ
(\sigma_{\Gamma_1^\dagger \cup \Sigma_1,\Gamma_2^\dagger} \cup \id_{\Sigma_2})
$$
and thus
\begin{multline*}
\ZZ(\m \cup \n)=(\id_{\ZZ(\Gamma_1) \tensor \ZZ(\Gamma_2)} \tensor
\ZZ(\widehat{\m \cup \n})) \circ (\ZZ(\ee_{\Gamma_1 \cup
\Gamma_2}) \tensor \id_{\ZZ(\Sigma_1)} \tensor
\id_{\ZZ(\Sigma_2)})\\
=\ZZ(\m) \tensor \ZZ(\n).
\end{multline*}
\item $\ZZ(\emptyset)=\FF.$  This is by definition.
\item $\ZZ(\sigma_{\Sigma,\Gamma})=\sigma_{\ZZ(\Sigma),\ZZ(\Gamma)}.$
This is also by definition.
\end{enumerate}
\end{proof}

\begin{rem}
From the point of view of the category theoretic definition of TQFT,
the distinction between sewing and mending may seem more semantic than
real.  In fact the mending property plays a distinct and crucial role.
Any $2$-framed or ordinary three-manifold invariant can be extended in a
formal fashion to something which satisfies all the other axioms (possibly at
the cost of having infinite-dimensional vector spaces). In fact Reshetikhin and
Turaev's original construction of the cobordism invariant  \cite{RT91}
demonstrated all the axioms but mending, relying on fairly general properties
of the link invariant and the modularity.  The demonstration of the mending
axiom
\cite{Turaev94,Sawin96a} uses properties of the link invariant which
have much more to do with the connection to conformal field theory.
Sewing is also the axiom which the Hennings invariant and its nonsemisimple
cousins fail.
\end{rem}

\section{Surgery}

The previous section reduced the problem of finding a TQFT from one
about $2$-framed cobordisms essentially to a question about $2$-framed
manifolds with boundary.  This is still quite difficult to work with,
and some combinatorial presentation is necessary.  The most convenient
is surgery.

First let us pick a convenient choice of representative objects 
$\Sigma_{g}.$  Specifically, let each $\Sigma_{g}$ be a 
handlebody of genus $g.$

Recall that an unoriented framed link in $S^3$ determines a (compact,
connected, simply-connected, oriented)
four-manifold with boundary as follows:  Identify $S^3$ with the
boundary of the four-ball $B^4,$ thicken the components of the link to
a collection of embedded tori in $S^3$ with a choice of longitude (the
framing), and attach a two-handle to $B^4$ along each component of the
link (the attaching part of the boundary of the two-handle has a
preferred longitude, so the attachment map is determined up to
isotopy).  If the original $S^3$ contains an embedding of handlebodies
and the link is chosen so as not to intersect the range of the
embeddings, an embedding of the same handlebodies is determined in the
boundary of the four-manifold.  Of course we may view the four-manifold 
as a $2$-framed three-manifold, and view the embedded handlebodies as 
its boundary.  Thus a framed link together with a collection of nonintersecting 
embeddings of handlebodies $\Sigma_g$ into the complement determines 
a basic $2$-framed cobordism.

\begin{figure}[hbt]
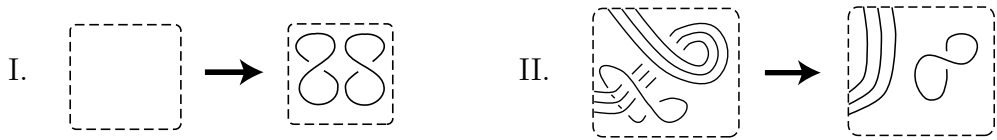


    $$\mbox{I.} \quad
\pic{fkirby1}{41}{-19}{0}{0}
\qquad \qquad
\mbox{II.} \quad
\pic{fkirby2}{50}{-22}{0}{0}$$ 
\caption{The framed Kirby moves with embedded handlebodies}\label{fg:fKirby}
\end{figure}

\begin{thm}
Every basic $2$-framed cobordism is diffeomorphic to one constructed by
surgery on $S^3$ with an embedding of handlebodies.  Two different
pairs of embeddings and surgery link in $S^3$ give the same $2$-framed 
cobordism if and only if they can be
connected by isotopy and a sequence of $2$-framed Kirby moves as shown
in Figure \ref{fg:fKirby}, where the open strands pictured in the
second move can be either strands of the surgery link or handles of
the embedded handlebodies (here the handlebodies are assumed to come
equipped with a choice of longitude for each handle, so that a
projection of an embedding can be presented by a projection of the
longitudes, a twist representing a Dehn twist as is the convention
with framed links).
\end{thm}

\begin{figure}[hbt]
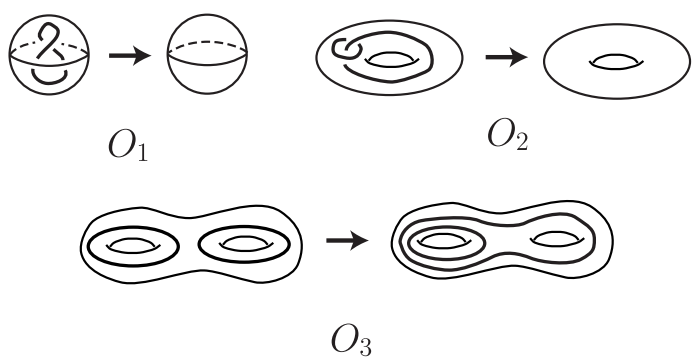

$$\pic{roberts}{130}{-60}{0}{0}$$
\caption{The Roberts moves}\label{fg:Roberts}
\end{figure}

\begin{figure}[hbt]
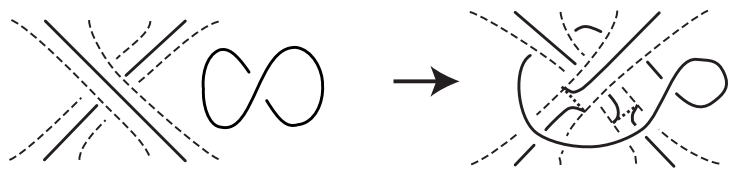

$$\pic{flip}{60}{-25}{0}{0}$$
\caption{Flipping a crossing with move II}\label{fg:flip}
\end{figure}

\begin{figure}[hbt]
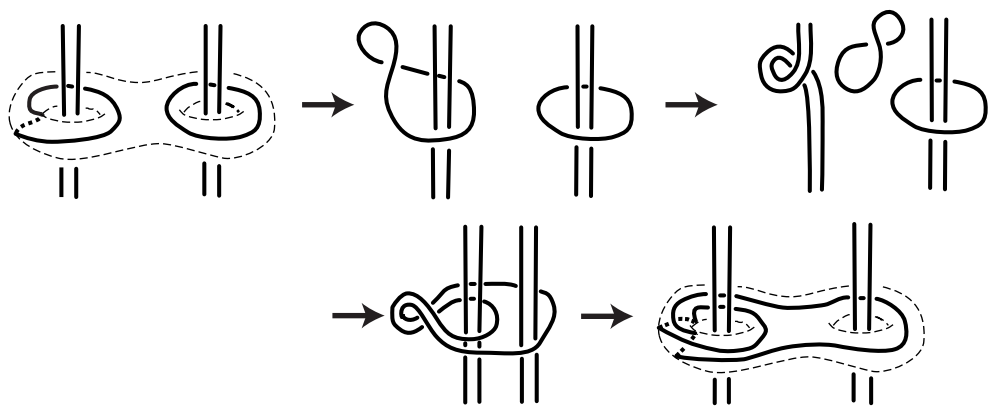

$$\pic{slide}{150}{-70}{0}{0}$$
\caption{Enacting $O_2$ with $2$-framed Kirby moves}\label{fg:handle}
\end{figure}

\begin{proof}
    
    Of course two basic $2$-framed cobordisms are the same if and 
    only if there is they have the same signature and there is a 
    diffeomorphism of the boundaries of the underlying four-manifold 
    intertwining the embeddings of each $\Sigma_g.$  
    
in \cite{Roberts97} Roberts proves that given two compact connected
three-manifolds with boundary and  a diffeomorphism from the boundary
of one to the boundary of the other, there is a framed link in the
interior of one such that surgery on that link gives a three-manifold
with boundary over which the diffeomorphism extends to a
diffeomorphism of the entire manifolds.  Further, two such links can
be connected by a sequence of moves $O_1,$ $O_2,$ $O_3$  shown in
Figure \ref{fg:Roberts} together with their mirror images and
inverses, where the (respectively) 
ball, two-handled torus and torus 
can be embedded anywhere in the manifold.   

We note first that in the proof of this result, move $O_1$ is used
only to make the signatures of the linking matrices of the two links
equal, after which every use of it can be a replaced with a use of move
$O_3,$ and thus the $2$-framed version of Roberts' result would simply
drop move $O_1.$  

Of course the situation in the present theorem is a special case of
Roberts' result, where one of the three-manifolds is $S^3$ with
several handlebodies removed, and the parameterizations of the
boundary  give the diffeomorphism.   Thus since it is clear that the
$2$-framed Kirby moves do not change the $2$-framed three-manifold with
embedded handlebodies, we need only prove that moves $O_2$ and $O_3$
can be replaced by a sequence of $2$-framed Kirby moves.

This relies on the following two observations, essentially due to Fenn
and Rourke \cite{FR79}.  The first is that,  given a link and
embedded handlebodies in $S^3,$ and given a choice of curves on the
boundary of the handlebodies which generate the homology of the
handlebody, a sequence of Kirby moves can replace the embedding with
one in which all these curves bound disks in $S^3$ minus the embedded
handlebodies.  The second is that if $K_1$ and $K_2$ are two links
plus embeddings in $S^3$ which can be connected by a sequence of
$2$-framed Kirby moves, and if $T$ is a framed link in one of the
embedded handlebodies, and $K_1',$ $K_2'$ are the result of embedding
$T$ in to $S^3$ via the embedding of the handlebody and then removing
that handlebody from the list, then $K_1'$ and $K_2'$ can be connected
by a sequence of $2$-framed Kirby moves.

For the first observation, notice the embedding of the handlebodies
gives an embedding of the boundary curves into $S^3,$ and as framed
links they admit a projection.   It is well-known that for any
projection of a framed link, one can, by flipping certain of the
crossings, make it a projection of a link in which all components are
unlinked 
zero framed unknots.  So apply $2$-framed Kirby move I to add
sufficiently many $\pm1$ framed unknot, and for each crossing that
needs to be flipped, apply move II as shown in Figure \ref{fg:flip}
to flip it, noting as indicated by the dotted lines in that figure
that the move is meant to be applied to the handlebodies, not just the
boundary curves.   The result of these moves will have the desired
property.

For the second observation, notice the result is  manifestly true
for a single $2$-framed Kirby move, and thus for an arbitrary sequence.

Figure \ref{fg:circum} decomposes move $O_3$ into a sequence of
$2$-framed Kirby moves assuming that the torus is embedded so that the
meridian plus the longitude bounds a disk in $S^3$ (here the vertical strands
represent any number of link components which might go through the 
bounded disks).
For an arbitrary embedding,  use the two observations to  show that a sequence
of
$2$-framed Kirby moves will  replace any embedding with one whose meridian plus
longitude bounds a disk, apply Figure
\ref{fg:circum}, then invert the sequence of Kirby moves to return to the
original embedding.

The two observations again reduce the general problem of move $O_2$ to
one where the two-handled torus is embedded as shown in Figure
\ref{fg:handle}, and the moves in the figure prove the result in that case.
\end{proof}

\begin{figure}[hbt]
$$\pic{circum}{70}{-32}{-30}{50}$$
\caption{Enacting $O_3$ with $2$-framed Kirby moves}\label{fg:circum}
\end{figure}

\begin{rem}
While \cite{Sawin96a} construct TQFTs by proving invariance of
appropriate quantities under the $2$-framed Kirby moves, one could just
as well use the Roberts moves directly.  
\end{rem}

\begin{figure}[hbt]
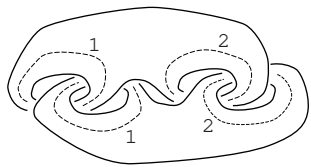

$$\pic{pairing}{60}{-25}{0}{0}$$
\caption{The embedding plus link corresponding to the pairing}\label{fg:pairing}
\end{figure}

\begin{figure}[hbt]
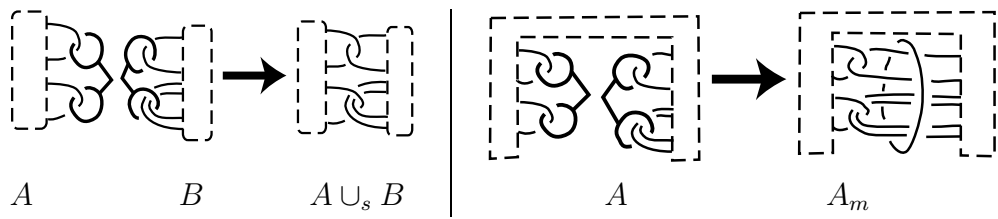
  \begin{tabular}{c|c}
$\begin{array}{c}\pic{sewing}{50}{-24}{0}{0}\\
\rule{0pt}{24pt} A \hspace{55pt} B \hspace{40pt}  A \cup_s B
\end{array}$ \hspace{20pt}&
$\begin{array}{c}
\pic{mending}{58}{-29}{0}{0}\\
\rule{0pt}{16pt}A \hspace{75pt} A_m
\end{array}$
\end{tabular}
\caption{Sewing and mending of three-manifolds}\label{fg:sewandmend}
\end{figure}

\begin{cor}
For each $g$ suppose $\Sigma_g$ is a handlebody of genus $g,$ in 
particular a $2$-framed surface, and suppose $\ZZ(\Sigma_g)$ is a finite-dimensional vector
space. Given
an embedding of  $\bigcup_{i=1}^n \Sigma_{g_i}$ with each
$\Sigma_{g_i}$ labeled by a vector in $\ZZ(\Sigma_{g_i})$ into $S^3$ together with a
framed unoriented link in the complement of the embedding, suppose $f$
is an invariant of the labeled embedding and the link.  Suppose
further that $f$ is linear as a function of each label, independent of
the ordering of $\{g_i\}_{i=1}^n,$ unchanged by the framed Kirby moves
and that the value of $f$ on the embedding shown in Figure
\ref{fg:pairing} (shown here for the genus two case) is a nondegenerate pairing
$\bracket{v,w}$ on
$\ZZ(\Sigma_g).$  Finally suppose that $f$ of the embeddings pictured
in Figure \ref{fg:sewandmend} are related by 
\begin{align*}
f(A \cup_s B)&= (f(A) \tensor f(B))\circ \Phi\\
f(A_m)&=f(A) \circ \Phi
\end{align*}
where $\Phi$ is constructed out of the dual element to the pairing as
in Theorem \ref{th:TQFT}.  Then $f$ is actually an invariant of the
basic cobordism determined by that embedding and link which satisfies the
axioms of Theorem \ref{th:TQFT} and as such determines a TQFT.
\end{cor}

\begin{proof}
Of course invariance under the biframed Kirby moves guarantees $f$ is an
invariant of the cobordism.  We are given that it satisfies the
Symmetry and Nondegeneracy Axioms of Theorem \ref{th:TQFT}, and to see
it satisfies Sewing and Mending it suffices to check that the
embeddings and link shown in Figure \ref{fg:sewandmend} (shown here
only for a genus two $2$-framed surface, represented by its underlying 
graph) represent the sewing and mending of the basic
cobordisms.  This follow easily from the definition of the surgery description.
\end{proof}

\bibliographystyle{alpha}

\end{document}